\documentclass[12pt]{amsart}

\usepackage{graphics, amsmath, amsfonts, amssymb, amscd, mathrsfs}

\newtheorem{theorem}{Theorem}
\newtheorem{proposition}[theorem]{Proposition}
\newtheorem{lemma}[theorem]{Lemma}

\theoremstyle{definition}
\newtheorem{definition}{Definition}
\newtheorem{example}[definition]{Example}

\newcommand{\A}{\mathscr{A}}
\newcommand{\B}{\mathscr{B}}
\newcommand{\C}{\mathbb{C}}
\newcommand{\D}{\mathbb{D}}

\newcommand{\R}{\mathbb{R}}

\newcommand{\T}{\mathbb{T}}
\newcommand{\PSH}{\operatorname{PSH}}
\renewcommand{\Re}{\operatorname{Re}}
\renewcommand{\Im}{\operatorname{Im}}

\title{Plurisubharmonic subextensions \\ as envelopes of disc functionals}

\author{Finnur L\'arusson}
\address{School of Mathematical Sciences, University of Adelaide, Adelaide SA 5005, Australia} 
\email{finnur.larusson@adelaide.edu.au}

\author{Evgeny A.\ Poletsky}
\address{Department of Mathematics, Syracuse University, 215 Carnegie Hall, Syracuse NY 13244, U.S.A.}
\email{eapolets@syr.edu}



\thanks{F.L.\ was supported by Australian Research Council grant DP120104110.  E.P.\ was supported by National Science Foundation grant DMS-0900877.  Work on this paper was begun when E.P.\ visited Adelaide in 2008 and continued during F.L.'s visit to Syracuse in 2012.  The authors thank each other's institutions for their hospitality.}

\date{27 January 2012.  Most recent changes 4 May 2012}

\begin{document}

\begin{abstract} 
We prove a disc formula for the largest plurisubharmonic subextension of an upper semicontinuous function on a domain $W$ in a Stein manifold to a larger domain $X$ under suitable conditions on $W$ and $X$.  We introduce a related equivalence relation on the space of analytic discs in $X$ with boundary in $W$.  The quotient, if it is Hausdorff, is a complex manifold with a local biholomorphism to $X$.  We use our disc formula to generalise Kiselman's minimum principle.  We show that his infimum function is an example of a plurisubharmonic subextension.
\end{abstract}

\maketitle
\tableofcontents

\section{Introduction}

\noindent
The theory of disc functionals was founded just over twenty years ago.  Its main goal is to provide disc formulas for important extremal plurisubharmonic functions in pluri\-potential theory, that is, to describe these functions as envelopes of disc functionals.  This brings the geometry of analytic discs into play in pluripotential theory.  Disc formulas have been proved for largest plurisubharmonic minorants (some of the main references, in historical order, are \cite{Poletsky1991}, \cite{BuSchachermayer1992}, \cite{Poletsky1993}, \cite{LarussonSigurdsson1998}, \cite{Rosay2003}, \cite{Magnusson2011}) and for various Green functions (see for example \cite{PoletskyShabat1989}, \cite{Edigarian1997}, \cite{LarussonSigurdsson2003}, \cite{RashkovskiiSigurdsson2005}, \cite{LarussonSigurdsson2009}).  For recent generalisations to singular spaces, see \cite{DrnovsekForstneric2011a}, \cite{DrnovsekForstneric2011b}.

We continue this project by proving a disc formula for largest plurisubharmonic subextensions.  Consider domains $W\subset X$ in complex affine space $\C^n$ or, more generally, in a Stein manifold.  Let $\phi:W\to[-\infty,\infty)$ be upper semicontinuous, for example plurisubharmonic, and let $S\phi$ be the supremum of all plurisubharmonic functions $u$ on $X$ with $u|W\leq\phi$.  If $X$ is covered by analytic discs with boundaries in $W$, then $S\phi$ is a plurisubharmonic function on $X$, the largest plurisubharmonic subextension of $\phi$ to $X$.  Under suitable conditions on $W$ and $X$, we prove that for every $x\in X$, $S\phi(x)$ is the infimum of the averages of $\phi$ over the boundaries of all analytic discs in $X$ with boundary in $W$ and centre $x$ (Theorems \ref{t:first-disc-formula} and \ref{t:second-disc-formula}).  In general, however, the disc formula can fail (Example \ref{e:counterexample}).

A recent Stein neighbourhood theorem of Forstneri\v c (\cite{Forstneric2007}, Theorem 1.2) allows us to work with analytic discs that are merely continuous on the closed unit disc $\overline\D$.  This is technically easier than the traditional approach that uses germs of holomorphic maps from open neighbourhoods of $\overline\D$.

A new equivalence relation on the space $\A_X^W$ of analytic discs in $X$ with boundary in $W$ naturally appears in the proof of our disc formula.  We call analytic discs in $\A_X^W$ {\it centre-homotopic} if they have the same centre and can be joined by a path in $\A_X^W$ of discs with that same centre.  The quotient of $\A_X^W$ by this equivalence relation, if it is Hausdorff, is a complex manifold with a local biholomorphism to $X$ (Theorem \ref{t:local-homeo}).  The sufficient conditions in Theorems \ref{t:first-disc-formula} and \ref{t:second-disc-formula} are naturally expressed in terms of the quotient (Theorem \ref{t:new-formulation}).  Finally, we use our disc formula to generalise Kiselman's minimum principle \cite{Kiselman1978} and to give a new proof of a special case of it, based on the observation that Kiselman's infimum function may be viewed as a plurisubharmonic subextension to a suitable larger domain (Theorem \ref{t:generalisation}).

\smallskip\noindent
{\it Acknowledgement.}  We thank Barbara Drinovec Drnov\v sek for pointing out an error in a previous version of the paper.

\section{A disc formula for plurisubharmonic subextensions}

\noindent
We start by establishing some notation.  For $r>0$, let $D_r=\{z\in\C:\lvert z\rvert<r\}$ and $\D=D_1$.  Let $\lambda$ denote the normalised arc length measure on the unit circle $\T=\partial\D$.  For a complex manifold $X$, let $\A_X$ denote the set of analytic discs in $X$, here taken to be continuous maps $f:\overline\D\to X$ that are holomorphic on $\D$.  We call $f(0)$ the centre of $f$.  We endow $\A_X$ with the compact-open topology, that is, the topology of uniform convergence on $\overline\D$.  It makes $\A_X$ a complete metrisable space.  For any topological space $Y$, a continuous map $Y\to \A_X$ is nothing but a continuous map $Y\times\overline\D\to X$ that is holomorphic when restricted to $\{y\}\times\D$ for every $y\in Y$.  If $W\subset X$, write $\A_X^W$ for the set of analytic discs $f$ in $X$ with $f(\T)\subset W$.  If $W$ is open, then $\A_X^W$ is open in $\A_X$. 

Let $f\in\A_X$.  By a recent theorem of Forstneri\v c (\cite{Forstneric2007}, Theorem 1.2), the graph $\Gamma_f=\{(z,f(z)):z\in\overline\D\}$ of $f$ has a basis of nice Stein open neighbourhoods in $\C\times X$.  More precisely, there is a basis of Stein open neighbourhoods $V$ of $\Gamma_f$ in $\C\times X$, each with a biholomorphism onto an open subset of $\C\times\C^{\dim X}$, mapping $(\{z\}\times X)\cap V$ onto an open convex subset of $\{z\}\times\C^{\dim X}$ for each $z\in\C$.  The sets $V^*=\{g\in\A_X: \Gamma_g\subset V\}$, as $V$ ranges over such a basis of open neighbourhoods of $\Gamma_f$, form a basis of open neighbourhoods of $f$ in $\A_X$.  It follows that there is a neighbourhood $W$ of $f(0)$ in $X$ and a continuous map $F:W\to\A_X$ such that $F(f(0))=f$ and $F(x)(0)=x$ for each $x\in W$.  Hence the centre map $c:\A_X\to X$, $f\mapsto f(0)$, is not only continuous but also open.  Each neighbourhood $V^*$ is contractible, so $\A_X$ is locally contractible.  In particular, the connected components and the path components of $\A_X$ are the same, and they are open in $\A_X$.

For an upper semicontinuous function $\phi:X\to[-\infty,\infty)$, let $H_\phi:\A_X\to[-\infty,\infty)$ be the {\it Poisson functional} associated to $\phi$, defined by the formula
\[ H_\phi(f) = \int_\T \phi\circ f\,d\lambda. \]
For $\B\subset\A_X$, the {\it Poisson envelope} $E_\B\phi:X\to[-\infty,\infty]$ of $\phi$ with respect to $\B$ is defined by the formula
\[ E_\B\phi(x)=\inf_{\begin{subarray}{c} f\in\B \\ f(0)=x \end{subarray}} H_\phi(f). \]
It is well known that $P\phi=E_{\A_X}\phi$ is the largest plurisubharmonic minorant of $\phi$ on $X$ (see the references given in the introduction).  If $W\subset X$ is open, $\B\subset\A_X^W$, and $\phi:W\to[-\infty,\infty)$ is upper semicontinuous, then the envelope $E_\B\phi:X\to[-\infty,\infty]$ is defined as above.

Let $W$ be a domain (a connected, nonempty, open subset) in a complex manifold $X$.
Let $\phi:W\to[-\infty,\infty)$ be upper semicontinuous, for example plurisubharmonic (we take the constant function $-\infty$ to be plurisubharmonic).  A plurisubharmonic function $u$ on $X$ is a {\it subextension} of $\phi$ if $u|W\leq\phi$.  Let 
\[ S\phi=\sup\{u\in\PSH(X):u|W\leq\phi\}. \]
Now $S\phi$ is upper semicontinuous (in particular nowhere $\infty$) and hence plurisubharmonic if and only if the class $\{u\in\PSH(X):u|W\leq\phi\}$ has local upper bounds on $X$, which holds for example if $X$ is covered by analytic discs with boundaries in $W$.  Then $S\phi$ is the largest plurisubharmonic subextension of $\phi$ to $X$.  It is easily seen that $S(S\phi)=S\phi$ and that $S\phi$ is maximal on $X\setminus\overline W$.  Also, $S\phi \leq E_{\A_X^W}\phi$, with equality if and only if $E_{\A_X^W}\phi$ is plurisubharmonic on $X$.  We will prove the disc formula $S\phi=E_{\A_X^W}\phi$ under suitable conditions on $X$ and $W$.

As a start, let us derive a preliminary disc formula for largest plurisubharmonic subextensions directly from the disc formula for largest plurisubharmonic minorants.

\begin{proposition}  \label{p:preliminary-disc-formula}
Let $W$ be a domain in a complex manifold $X$, such that $X$ is covered by analytic discs with boundaries in $W$.  Let $\phi:W\to[-\infty,\infty)$ be upper semicontinuous and bounded above.  Then, for every $x\in X$,
\[ S\phi(x)=\lim_{\epsilon\to 0+}\inf_f \int_{(f|\T)^{-1}(W)} \phi\circ f\,d\lambda, \]
where, for each $\epsilon>0$, the infimum is taken over all $f\in\A_X$ such that $f(0)=x$ and $\lambda((f|\T)^{-1}(W))>1-\epsilon$.
\end{proposition}

\begin{proof}
Suppose $\phi$ is bounded above on $W$ by $m\in\mathbb N$.  For each $n\geq m$, define an upper semicontinuous function $\phi_n$ on $X$ as $\phi$ on $W$ and as $n$ on $X\setminus W$.  Now $P\phi_n$ is plurisubharmonic on $X$ and $P\phi_n\leq \phi$ on $W$, so $P\phi_n\leq S\phi$.  Also, $S\phi\leq m$ on $X$, so $S\phi\leq P\phi_n$.  Hence $S\phi=P\phi_n$ for all $n\geq m$.

Let $x\in X$.  If $f\in\A_X$ has $\lambda((f|\T)^{-1}(W))>1-\epsilon$ and $f(0)=x$, then
\[ S\phi(x)\leq \int_\T\phi_m\circ f\,d\lambda < \int_{(f|\T)^{-1}(W)} \phi\circ f\,d\lambda + m\epsilon. \]
Thus $S\phi(x)\leq \lim\limits_{\epsilon\to 0+}\inf\limits_f \int_{(f|\T)^{-1}(W)} \phi\circ f\,d\lambda$.  On the other hand, for each $n\geq m$, there is $f_n\in\A_X$ with $f_n(0)=x$ and $H_{\phi_n}(f_n)\leq P\phi_n(x)+\tfrac 1 n = S\phi(x)+\tfrac 1 n$.  Then $\lambda((f_n|\T)^{-1}(W))\to 1$ and
\[ \int_{(f_n|\T)^{-1}(W)} \phi\circ f_n\,d\lambda \leq \int_\T \phi_n\circ f_n\,d\lambda \leq S\phi(x)+\tfrac 1 n. \]
Thus $S\phi(x)\geq \lim\limits_{\epsilon\to 0+}\inf\limits_f \int_{(f|\T)^{-1}(W)} \phi\circ f\,d\lambda$.
\end{proof}

The disc formula in Proposition \ref{p:preliminary-disc-formula} is rather clumsy.  It is natural to ask whether we can \lq\lq pass to the limit\rq\rq\ and use only analytic discs whose entire boundary lies in $W$.  This turns out to be a subtle question, involving a mix of complex analysis and topology.  The answer is affirmative only if suitable restrictions are imposed on the pair $W\subset X$.

We say that $f_0, f_1\in\A_X^W$ with $f_0(0)=f_1(0)$ are {\it centre-homotopic} if there is a continuous map $f:\overline \D\times [0,1]\to X$ with $f(\cdot,t)\in\A_X^W$ for all $t\in [0,1]$ and $f(\cdot,t)=f_t$ for $t=0,1$ (that is, a continuous path in $\A_X^W$ joining $f_0$ and $f_1$), such that $f(0,t)=f_0(0)$ for all $t\in[0,1]$.  

A $W$-{\it disc structure} on $X$ is a family $\beta=(\beta_\nu)$ of continuous maps $\beta_\nu:U_\nu\to \A_X^W$, where $(U_\nu)$ is an open cover of $X$, satisfying the following two conditions.
\begin{itemize}
\item[]  For all $x\in U_\nu$, $\beta_\nu(x)(0)=x$.
\item[(S)]  If $x\in U_\nu\cap U_\mu$, then $\beta_\nu(x)$ and $\beta_\mu(x)$ are centre-homotopic.
\end{itemize}
We will be interested in the following condition that $\beta$ may or may not satisfy.
\begin{itemize}
\item[(N)]  There is $\mu$ such that $U_\mu=W$ and $\beta_\mu(w)$ is the constant disc at $w$ for all $w\in W$.
\end{itemize}
The class $\B_\beta\subset\A_X^W$ associated to a $W$-disc structure $\beta$ on $X$ is the union $\bigcup\limits_\nu\beta_\nu(U_\nu)$.  It is easily seen that if $\phi$ is upper semicontinuous on $W$, then the envelope $E_\beta\phi=E_{\B_\beta}\phi$ is upper semicontinuous on $X$.  We say that $X$ is a {\it schlicht disc extension} of $W$ if $X$ carries a $W$-disc structure satisfying (N).  These definitions will be viewed in a more abstract light in the following section.

\begin{example}  \label{e:simple-example}
Here is a simple example to illustrate the definitions.  For $n\geq 2$ and $r>0$, let $B_r^n=\{x\in\C^n:\lVert x\rVert<r\}$.  Set $X=B_4^n$ and $W=B_4^n\setminus\overline {B_1^n}$.  Let $U_0=W$ and $U_1=B_2^n$.  Then $\{U_0, U_1\}$ is an open cover of $X$, and a $W$-disc structure on $X$ satisfying (N) is given by setting $\beta_0(x)(\zeta)=x$ for all $x\in U_0$, and
\[ \beta_1(x)(\zeta) =\big(\rho(x)\dfrac{\rho(x)\zeta+x_1}{\rho(x)+\bar x_1\zeta},x_2,\dots,x_n\big), \]
where
\[ \rho(x)=\sqrt{9-\lvert x_2\rvert^2-\cdots-\lvert x_n\rvert^2}, \]
for all $x=(x_1,\dots,x_n)\in U_1$.  Condition (S) is evident.
\end{example}

The following lemma is proved along the lines of the original proof in \cite{Poletsky1991} of the plurisubharmonicity of the Poisson envelope.  We follow the exposition in \cite{LarussonSigurdsson1998}.  Note that we are now restricting our discussion to domains in affine space.

\begin{lemma}  \label{l:inequality}
Let $W\subset X$ be domains in $\C^n$, and $\beta$ be a $W$-disc structure on $X$.  If $\phi:W\to[-\infty,\infty)$ is upper semicontinuous, then
\[ E_{\A_X^W}\phi \leq P E_\beta \phi. \]
\end{lemma}

Before proving the lemma we state and prove our first theorem.  

\begin{theorem}  \label{t:first-disc-formula}
Let $W\subset X$ be domains in $\C^n$ such that $X$ is a schlicht disc extension of $W$.  If $\phi:W\to[-\infty,\infty)$ is upper semicontinuous, then
\[ S\phi = E_{\A_X^W}\phi. \]
\end{theorem}

\begin{proof}
Let $\beta$ be a $W$-disc structure on $X$ satisfying (N).  We have
\[ S\phi \leq E_{\A_X^W}\phi \leq P E_\beta \phi\leq S\phi. \]
The first inequality is obvious.  The second is Lemma \ref{l:inequality}.  The third holds because $E_\beta\phi\leq\phi$ on $W$ by (N).
\end{proof}

This approach to proving a disc formula using an auxiliary class such as $\B_\beta$ first appeared in \cite{LarussonSigurdsson2005}.

\begin{proof}[Proof of Lemma \ref{l:inequality}]
Let $\beta$ be a $W$-disc structure on $X$.  To prove the desired inequality, we show that for every $h\in\A_X$, $\epsilon>0$, and a continuous function $v:X\to\R$ with $v\geq E_\beta\phi$, there exists $g\in\A_X^W$ such that $g(0)=h(0)$ and 
\begin{equation}\label{ineq:1} H_\phi(g)\leq H_v(h)+\epsilon.  \end{equation}
The proof is carried out in three steps.  First we show that there exists a continuous map $F:\overline\D\times \T\to X$, such that $F(\cdot,w)\in \A_X^W$ and $F(0,w)=h(w)$ for all $w\in \T$, and
\begin{equation}\label{ineq:2} \int_\T H_\phi(F(\cdot, w))\, d\lambda(w)\leq H_v(h)+\epsilon/2. \end{equation}
Next we show that there exists a continuous map $G:\overline\D\times \overline\D\to X$, holomorphic on $\D\times\D$, such that $G(0,w)=h(w)$  for all $w\in\overline\D$, $G(\T\times\T)\subset W$, and
\begin{equation}\label{ineq:3} \int_\T H_\phi(G(\cdot, w))\, d\lambda(w) \leq  \int_\T H_\phi(F(\cdot, w))\, d\lambda(w) +\epsilon/2. \end{equation}
Finally we show that there is $\theta_0\in[0,2{\pi}]$ such that if $g\in\A_X^W$ is defined by the formula $g(z)=G(e^{i\theta_0}z,z)$, then
\begin{equation}\label{ineq:4}  H_\phi(g)\leq \int_\T H_\phi(G(\cdot, w))\, d\lambda(w). \end{equation}
By combining (\ref{ineq:2}), (\ref{ineq:3}), and (\ref{ineq:4}), we get (\ref{ineq:1}).

\smallskip
\noindent
{\it Step 1.}  Let $h\in\A_X$, $\epsilon>0$, and $v:X\to\R$ be continuous with $v\geq E_\beta\phi$.  Let $w_0\in \T$ and set $x_0=h(w_0)$.  Find $\nu$ such that $x_0\in U_\nu$ and $H_\phi(\beta_\nu(x_0))<v(x_0)+\epsilon/(8\pi)$.  For all $x$ in a small enough neighbourhood of $x_0$,  we have $H_\phi(\beta_\nu(x))<v(x)+\epsilon/(8\pi)$.  By compactness, there is a cover of $\T$ by closed arcs $I_1,\dots,I_m$ that meet only in end points, such that for each $j=1,\dots,m$, there is $\nu_j$ with $h(I_j)\subset U_{\nu_j}$ and $H_\phi(\beta_{\nu_j}(h(w)))<v(h(w))+\epsilon/(8\pi)$ for all  $w\in I_j$.

Let $w_0$ be a common end point of $I_1$ and $I_2$, say.  By (S), $\beta_{\nu_1}(h(w_0))$ and $\beta_{\nu_2}(h(w_0))$ are centre-homotopic.  Choose a centre-homotopy between them, use small translations to spread the centres of the analytic discs in the homotopy over a small arc centred at $w_0$, and reparametrise $\beta_{\nu_1}\circ h$ and $\beta_{\nu_2}\circ h$ by small translations over subarcs slightly smaller than $I_1$ and $I_2$, respectively.  In this way we obtain a continuous map $F:\overline\D\times \T\to X$, such that $F(\cdot,w)\in \A_X^W$ and $F(0,w)=h(w)$ for each $w\in\T$, and
\begin{align*}
\int_\T H_\phi(F(\cdot, w))\, d\lambda(w) &\leq \sum_{j=1}^m  \int_{ I_j} H_\phi(\beta_{\nu_j}\circ h)\, d\lambda+\epsilon/4  \\ &\leq \sum_{j=1}^m \int_{I_j} v\circ h\, d\lambda +\epsilon/2 =  H_v(h)+ \epsilon/2,
\end{align*}
so we have proved (\ref{ineq:2}).

\smallskip
\noindent
{\it Step 2.}  For each $j\geq 1$, we define the Ces\`aro mean $F_j:\overline\D\times \overline\D^*\to\C^n$, where $\overline\D^*$ denotes $\overline\D\setminus\{0\}$, by
\[ F_j(z,w)= h(w)+\dfrac 1{j+1}\sum_{m=0}^j\sum_{k=-m}^m \bigg(\int_\T \big(F(z, \zeta)-h(\zeta)\big) \zeta^{-k}\, d\lambda(\zeta) \bigg) w^k. \]
The well-known theorem on the uniform convergence of the Ces\`aro means of a continuous function on $\T$ holds for maps into a Banach space such as the space of continuous maps $\overline\D\to\C^n$ that are holomorphic on $\D$.  We conclude that $F_j\to F$ uniformly on $\overline\D\times \T$ as $j\to \infty$.  Hence $F_j(\overline\D\times\T)\subset X$ and $F_j(\T\times\T)\subset W$ for $j$ large enough.  For simplicity, assume that this holds for all $j\geq 1$.

For every $z\in\overline\D$, the map $w\mapsto F_j(z,w)-h(w)$ has a pole of order at most $j$ at the origin, and for every $w\in \overline\D^*$, the map $z\mapsto F_j(z,w)-h(w)$ has a zero at the origin.  Thus $(z,w)\mapsto F_j(zw^k,w)$ extends to a continuous map $\overline \D\times \overline \D \to \C^n$, holomorphic on $\D\times\D$, for every $k\geq j$.  

Since $F_j(0,w)=h(w)\in X$ for all $w\in\overline\D^*$, there is $\delta_j >0$ such that $F_j(zw^k,w)\in X$ for all integers $k\geq j$ and  $(z,w)\in D_{\delta_j}\times \overline \D$.  Since $F_j(\overline \D\times \T)\subset X$, there is $r_j\in (0,1)$ such that $F_j(\overline \D\times (\overline \D\setminus D_{r_j}))\subset X$, so $F_j(zw^k,w)\in X$ for all $(z,w)\in \overline \D\times (\overline \D\setminus D_{r_j})$ and all $k\geq j$.    

Take $k_j\geq j$ so large that $\lvert zw^{k_j}\rvert<\delta_j$ for all $(z,w)\in \overline \D\times D_{r_j}$.  Then $F_j(zw^{k_j},w)\in X$ for all $(z,w)\in\overline \D\times \overline \D$.  Define a continuous map $G_j:\overline\D\times \overline\D\to X$, holomorphic on $\D\times\D$, by $G_j(z,w)=F_j(zw^{k_j},w)$.  Then $G_j(\T\times\T)\subset W$ and $G_j(0,w)=h(w)$ for all $w\in \overline\D$.

Take $j$ large enough that
\begin{align*} \dfrac 1{(2\pi)^2}\int_0^{2{\pi}}\int_0^{2\pi} \phi(F_j(e^{it},e^{i\theta}))\, dtd{\theta} &\leq \dfrac 1{(2\pi)^2}\int_0^{2{\pi}}\int_0^{2{\pi}} \phi(F(e^{it},e^{i{\theta}}))\, dtd{\theta}
+\epsilon/2 \\  &=\int_0^{2{\pi}} H_\phi(F(\cdot, w))\, d\lambda(w)+\epsilon/2,
\end{align*}
and let $G=G_j$.  Then
\begin{align*} \int_\T H_\phi(G(\cdot, w))\, d\lambda(w) &=  \dfrac 1{(2\pi)^2}\int_0^{2\pi}\int_0^{2\pi} \phi(F_j(e^{i(t+k_j{\theta})},e^{i{\theta}}))\, dtd{\theta} \\ &=  \dfrac 1{(2\pi)^2}\int_0^{2\pi}\int_0^{2\pi} \phi(F_j(e^{it},e^{i\theta}))\, dtd\theta,
\end{align*}
and (\ref{ineq:3}) is proved.

\smallskip
\noindent
{\it Step 3.}  The right-hand side of (\ref{ineq:4}) is
\[ \dfrac 1{(2\pi)^2}\int_0^{2{\pi}}\int_0^{2{\pi}} \phi(G(e^{it},e^{i\theta}))\, dt\, d\theta = \dfrac 1{(2\pi)^2}\int_0^{2{\pi}}\int_0^{2{\pi}} \phi(G(e^{i\theta}e^{it},e^{it}))\, dt\, d\theta. \]
There is $\theta_0\in [0,2\pi]$ such that
\[ \dfrac 1{(2\pi)^2}\int_0^{2{\pi}}\int_0^{2\pi} \phi(G(e^{it},e^{i\theta}))\, dtd\theta \geq \dfrac 1{2\pi}\int_0^{2\pi} \phi(G(e^{i\theta_0}e^{it},e^{it})) \, dt. \]
If we set $g(z)=G(e^{i\theta_0}z,z)$, then $g(0)=G(0,0)$ and (\ref{ineq:4}) holds.
\end{proof}

We now provide another sufficient condition for our disc formula to hold.

\begin{theorem}  \label{t:second-disc-formula}
Let $W\subset X$ be domains in $\C^n$.  Suppose $\A_X^W$ has a connected component, call it $\B$, with the following properties.
\begin{enumerate}
\item  $\B$ covers $X$.
\item  If two analytic discs in $\B$ have the same centre, then they are centre-homotopic.
\end{enumerate}
Then, for every upper semicontinuous function $\phi:W\to[-\infty,\infty)$,
\[ S\phi = E_{\A_X^W}\phi. \]
\end{theorem}

\begin{proof}
The proof is the same as the proof of Theorem \ref{t:first-disc-formula}, except for the inequality $P E_\B \phi\leq S\phi$, which now requires more work.  We can show that $E_{\A_X^W}\phi \leq P E_\B \phi$ exactly as we proved Lemma \ref{l:inequality}, using (1) and (2) and the fact that $\B$ is open in $\A_X^W$.  We need to show that $E_\B\phi\leq\phi$ on $W$, which previously was an immediate consequence of (N).

Let $p\in W$ and $\epsilon>0$.  We need an analytic disc $f\in\B$ with $f(0)=p$ and $H_\phi(f)<\phi(p)+\epsilon$.  Take $g\in\B$ with $g(0)=p$.  Extend $g$ to a continuous map $g:\overline\D\cup[1,2]\to X$ such that $g|[1,2]$ is a path in $W$ with $g(2)=p$.  By Mergelyan's theorem, $g$ can be approximated uniformly on $\overline\D\cup[1,2]$ by polynomial maps $\C\to\C^n$.  Since $W$ and $\B$ are open, we may assume that $g$ is the restriction to $\overline\D\cup[1,2]$ of a polynomial map that we will still call $g$.

Let $\Omega\subset\C$ be the simply connected domain of all points within distance $\delta>0$ of $\overline\D\cup[1,2]$.  Let $\mu$ be the harmonic measure of $\Omega$ with respect to the point $2$.  Choose $\delta$ so small that $g(\overline\Omega\setminus\D)\subset W$ and $\int_{\partial\Omega}\phi\circ g\,d\mu<\phi(g(2))+\epsilon=\phi(p)+\epsilon$.  

A theorem of Rad\'o (\cite{Rado1922}; \cite{Goluzin1969}, II.5, Theorem 2) states that as a simply connected bounded domain in $\C$ is continuously varied in a suitable sense, its normalised Riemann map also varies continuously.  More precisely, for each $n\geq 0$, let $U_n$ be a simply connected domain in $\C$ containing $0$ and bounded by Jordan curves, and let $\psi_n:\D\to U_n$ be the biholomorphism with $\psi_n(0)=0$ and $\psi_n'(0)>0$.  Then $\psi_n\to \psi_0$ uniformly on $\overline\D$ if and only if for every $\epsilon>0$, there is $m\geq 1$, such that for every $n\geq m$, there is a homeomorphism $\partial U_n\to\partial U_0$ such that the distance between corresponding points is at most $\epsilon$.

Rad\'o's theorem provides a continuous map $\Psi:\overline\D\times[0,1]\to\overline\Omega$ such that $\Psi_t=\Psi(\cdot,t)$ is a homeomorphism onto its image and holomorphic on $\D$ with $\Psi_t(0)=0$ for every $t\in[0,1]$, $\Psi(\T\times[0,1])\subset\overline\Omega\setminus\D$, $\Psi_0$ is the inclusion $\overline\D\hookrightarrow\overline\Omega$, and $\Psi_1(\overline\D)=\overline\Omega$.  Then $g\circ\Psi_t\in\A_X^W$ for all $t\in[0,1]$, so $h=g\circ\Psi_1\in\B$ and $h(0)=p$.  Let $a=\Psi_1^{-1}(2)$.  Then $h(a)=p$.  By precomposing $h$ by a path in $\operatorname{Aut}\D$ joining the identity to an automorphism $\alpha$ of $\D$ taking $0$ to $a$, we obtain $f\in\B$ with $f(0)=p$.  Finally, 
\[ H_\phi(f)=\int_\T\phi\circ g\circ\Psi_1\circ\alpha\,d\lambda=\int_{\partial\Omega}\phi\circ g\,d\mu<\phi(p)+\epsilon. \qedhere \]
\end{proof}

In Theorems \ref{t:first-disc-formula} and \ref{t:second-disc-formula}, $\C^n$ can be replaced by any Stein manifold $S$.  Only minor modifications of Step 2 in the proof of Lemma \ref{l:inequality} and of the proof of Theorem \ref{t:second-disc-formula} are needed, using a tubular neighbourhood of $S$ viewed as a submanifold of $\C^m$ for some $m$.  Whether $\C^n$ can be replaced by an arbitrary complex manifold is an open question.

\begin{example}
Take $X=D_2$ and $W=D_2\setminus\overline\D$.  Then we have one connected component of $\A_X^W$ for each nonnegative winding number.  The hypotheses of Theorem \ref{t:second-disc-formula} hold for each positive winding number.  Namely, (1) is obvious and (2) follows from writing an analytic disc in $\A_X^W$ as the product of a holomorphic function without zeros and a Blaschke product whose degree equals the winding number.  Thus our disc formula holds for the pair $W\subset X$.  Theorem \ref{t:first-disc-formula} does not apply since $X$ is not a schlicht disc extension of $W$.
\end{example}

The following example shows that without the hypotheses of Theorems \ref{t:first-disc-formula} and \ref{t:second-disc-formula}, our disc formula can fail.  

\begin{example}  \label{e:counterexample}
Fix $\delta\in(0,\tfrac 1 2)$ and let
\begin{align*} W_1&=\D\times\{z_2\in\C: \lvert z_2\rvert<\delta\},\\
W_2&=\D\times\{z_2\in\C: 1-\delta<\lvert z_2\rvert<1\}. \end{align*}
Join $W_1$ and $W_2$ by the curve $[0,1]\to\C^2$, $t\mapsto (1+e^{2\pi i(2t-1)/3}, (1-\frac \delta 2)t)$.  Let $W_3$ be a thin open tubular neighbourhood of the image of the curve, and let $W$ be the domain $W_1\cup W_2\cup W_3$.  We may assume that the intersection with $\R$ of the projection of $W$ onto the $z_1$-plane is $(-1,1)\cup (a,b)$ with $1<a<2<b$, and that $b$ is the supremum of $\Re z_1$ on $W$.  Let $I=(a,b)$,
\begin{align*} W^-&=W_1\cup\{z\in W_3:\Im z_1<0\}, \\
W^+&=W_2\cup\{z\in W_3:\Im z_1>0\}. \end{align*}
Let $U_1$ and $U_2$ be thin open neighbourhoods of the semicircles
\[ \gamma_1=\{e^{it}/2:t\in[0,\pi]\} \quad\textrm{and}\quad \gamma_2=\{e^{it}/2:t\in[\pi,2\pi]\} \]
in $\D$, respectively, and let
\begin{align*} V_1&=U_1\times \{z_2\in\C: \lvert z_2\rvert<\delta\}\subset W_1, \\
V_2&=U_2\times\{z_2\in\C: 1-\delta<\lvert z_2\rvert<1\}\subset W_2. \end{align*}
Define an upper semicontinuous function $\phi$ on $W$ as $-1$ on $V_1\cup V_2$, and as $0$ elsewhere.

Let $X$ be any domain in $\C^2$ containing $W$, such that $X$ is covered by images $f(\overline\D)$ of analytic discs $f\in\A_{\C^2}^W$ and such that $X$ contains the image of the analytic disc $\zeta\mapsto(z_1,(1-\tfrac \delta 2)\zeta)$ for every $z_1\in \gamma_2$.  For example, $X$ could be $\D^2\cup W_3$.

By making $W_3$, $U_1$, and $U_2$ thin enough, we can insure that there is $s<1$ such that for every $g\in\A_{\C^2}$ with $g(0)=(0,0)$, each of the sets $g^{-1}(I\cup V_1)\cap\T$ and $g^{-1}(I\cup V_2)\cap\T$ has harmonic measure (with respect to $0\in\D$, that is, normalised arc length measure) less than $s$.

We will show that $S\phi(0,0)=-1$, but $E_{\A_X^W}\phi(0,0)\geq -s >-1$, so our disc formula fails for the pair $W\subset X$.

First, if $z_1\in\gamma_1$, then $\phi(z_1,0)=-1$, so $S\phi(z_1,0)=-1$.  If $z_1\in \gamma_2$, then the analytic disc $\zeta\mapsto(z_1,(1-\tfrac \delta 2)\zeta)$ has its boundary in $V_2$, so $\phi=-1$ on the boundary, and $S\phi(z_1,0)=-1$.  It follows that $S\phi(0,0)=-1$.

Second, let $f=(f_1,f_2)\in\A_X^W$ with $f(0)=(0,0)$.  We may assume that $f$ extends holomorphically to an open neighbourhood of $\overline\D$, that $f_1$ has no critical values in $\bar I$, and that $f_1|\T$ is transverse to $I$, so in particular $f_1^{-1}(I)\cap\T$ is finite.  We claim that the harmonic measure of $f^{-1}(V_1\cup V_2)\cap\T$ is less than $s$, so $H_\phi(f)>-s$.

If $f_1^{-1}(I)=\varnothing$, then $f^{-1}(V_1\cup V_2)\cap\T$ is either $f^{-1}(V_1)\cap\T$ or $f^{-1}(V_2)\cap\T$, so the claim is clear.  Assume that $f_1^{-1}(I)\neq\varnothing$.  Then $f_1^{-1}(I)$ is the disjoint union of finitely many embedded $1$-dimensional submanifolds.  None of them are loops, for otherwise $f_1$ would be constant, so they are all arcs.  None of them are relatively compact in $\D$ since $f_1$ cannot take a point in $\D$ to $b$, so each arc has one end point on $\T$ and the other on $\T$ or in $\D$.  

Say an arc in $f_1^{-1}(I)$ is {\it good} if both its end points lie on $\T$.  Call the other arcs {\it bad}.  Let $\Omega$ be the set of all points in $\D$ that lie on the same side of each good arc as $0$.  Then $\Omega$ is a simply connected domain containing $0$, bounded by some of the good arcs and some arcs in $\T$ that we will call {\it circular arcs}, unless there are no good arcs, in which case $\Omega=\D$ and $\T$ is the one circular \lq\lq arc\rq\rq.

The finitely many points on each circular arc that lie in $f_1^{-1}(I)$ divide the circular arc into open {\it subarcs}, each subarc having its $f$-image in $W^-$ or $W^+$.  Suppose there is a subarc with its $f$-image in $W^+$.  Moving counterclockwise along the subarc, we come to an end point $p$.  

Suppose that on the other side of $p$ there is a subarc of the same circular arc.  Then $p$ is an end point of a bad arc, and $\Im f_1|\T$ changes sign at $p$ from positive to negative.  Since $f_1$ is conformal at $p$, $\Re f_1$ is decreasing as we approach $p$ along the bad arc.  Now $\Re f_1$ has no critical points on the bad arc, so $\Re f_1$ is decreasing all along the bad arc in the direction of $p$, and the maximum principle is violated at the interior end point of the bad arc.  A slight elaboration of this argument shows that $\Omega\neq\D$.

Therefore $p$ is an end point of a good arc.  Since $f_1$ is conformal at $p$, $\Re f_1$ is increasing as we leave $p$ along the good arc.  Since $\Re f_1$ has no critical points on the good arc, $\Re f_1$ is still increasing as we approach the other end point $q$ of the good arc.  Since $f_1$ is conformal at $q$, $\Im f_1$ is increasing as we leave $q$ along the next subarc, which therefore has its $f$-image in $W^+$.

This shows that $f(\partial\Omega)$ cannot intersect both $W^-$ and $W^+$.  Thus either $f^{-1}(V_1)\cap\T$ or $f^{-1}(V_2)\cap\T$ lies behind $\partial\Omega\cap\D$ as seen from $0$, so the harmonic measure of $f^{-1}(V_1\cup V_2)\cap\T$ is at most the harmonic measure with respect to $0\in\Omega$ of either $f^{-1}(I\cup V_1)\cap\partial\Omega$ or $f^{-1}(I\cup V_2)\cap\partial\Omega$, and is therefore less than $s$, proving our claim.

By Theorem \ref{t:first-disc-formula}, $X$ is not a schlicht disc extension of $W$.  By Theorem \ref{t:second-disc-formula}, one or both of the following hold.
\begin{itemize}
\item  No disc in $\A_X^W$ can be deformed in $\A_X^W$ to a disc with an arbitrary centre in $X$.
\item  There are discs in $\A_X^W$ with the same centre that are homotopic but not centre-homotopic.
\end{itemize}
\end{example}

\section{The centre-homotopy relation on spaces of analytic discs}

\noindent
Let $X$ be a complex manifold and, as before, let $\A_X$ be the space of analytic discs in $X$, that is, continuous maps $\overline\D\to X$ that are holomorphic on $\D$, with the topology of uniform convergence.  Let $f\in\A_X$.  Recall that by (\cite{Forstneric2007}, Theorem 1.2), there is a basis of Stein open neighbourhoods $V$ of the graph $\Gamma_f$ of $f$ in $\C\times X$, each with a biholomorphism onto an open subset of $\C\times\C^{\dim X}$, mapping $(\{z\}\times X)\cap V$ onto an open convex subset of $\{z\}\times\C^{\dim X}$ for each $z\in\C$.  The sets $V^*=\{g\in\A_X: \Gamma_g\subset V\}$, as $V$ ranges over such a basis of open neighbourhoods of $\Gamma_f$, form a basis of open neighbourhoods of $f$ in $\A_X$.  We already noted that $V^*$ is contractible; moreover, $V^*$ intersects each fibre of the centre map $\A_X\to X$, $f\mapsto f(0)$, in a contractible set.  Let $W$ be a domain in $X$.  If $f\in\A_X^W$ and $V$ as above is small enough, then $V^*$ lies in $\A_X^W$ and intersects each fibre of the restricted centre map $c:\A_X^W\to X$ in a contractible set.

Define an equivalence relation $\sim_c$ on $\A_X^W$ by taking $f~\sim_c~g$ if $f$ and $g$ are centre-homotopic, that is, the equivalence classes of $\sim_c$ are the connected components, or equivalently the path components, of the fibres of $c$.  Let $q:\A_X^W\to X_W=\A_X^W/\sim_c$ be the quotient map, and endow $X_W$ with the quotient topology.  Let $X_W^0$ be the connected component of $X_W$ containing the equivalence classes of the constant discs in $W$.  Since $c$ is continuous, it factors through $q$ by a continuous map $\pi:X_W\to X$.

\begin{proposition}  \label{p:open-relation}
Let $W$ be a domain in a complex manifold $X$.  The equivalence relation $\sim_c$ on $\A_X^W$ is open, that is, the quotient map $q:\A_X^W\to X_W$ is open.
\end{proposition}

\begin{proof}
Let $f\sim_c g$ in $\A_X^W$ and let $V$ be a neighbourhood of $f$.  Take a path $\gamma :[0,1]\to c^{-1}(f(0))$ with $\gamma(0)=f$, $\gamma(1)=g$.  Cover the image of $\gamma$ by a chain of open sets $V_1^*,\dots,V_k^*$ of the kind described above, whose intersection with each fibre of $c$ is contractible, such that $f\in V_1^*\subset V$, $g\in V_k^*$, and $V_j^*\cap V_{j+1}^*\cap\gamma([0,1])\neq\varnothing$ for $j=1,\dots,k-1$. Recall that $c$ is open and let $U$ be the neighbourhood $\bigcap\limits_{j=1}^{k-1} c(V_j^*\cap V_{j+1}^*)$ of $f(0)$ in $X$.  Then each analytic disc in the neighbourhood $V_k^*\cap c^{-1}(U)$ of $g$ is $\sim_c$-related to an analytic disc in $V$.  Thus $\sim_c$ is open.
\end{proof}

\begin{theorem}  \label{t:local-homeo}
Let $W$ be a domain in a complex manifold $X$.  Then the map $\pi:X_W\to X$ is a local homeomorphism.
\end{theorem}

\begin{proof}
First, $\pi$ is open since $c$ is.  We need to show that $\pi$ is locally injective.  Let $f\in\A_X^W$.  Let $U$ be a neighbourhood of $q(f)$ in $X_W$.  Then $q^{-1}(U)$ is a neighbourhood of $f$ in $\A_X^W$.  Find a neighbourhood $V^*$ of $f$ in $\A_X^W$ as described above with $V^*\subset q^{-1}(U)$.  By Proposition \ref{p:open-relation}, $q(V^*)\subset U$ is a neighbourhood of $q(f)$, and $\pi$ is injective on $q(V^*)$.
\end{proof}

The Poincar\'e-Volterra theorem now implies that each connected component of $X_W$ is second countable (\cite{Bourbaki1966}, I.11.7, Corollary 2).  Thus we could turn $X_W$ into a complex manifold with the unique complex structure that makes $\pi$ holomorphic, except that we do not know whether $X_W$ is Hausdorff.  It is evident, though, that $X_W$ has closed points, that is, is $T_1$.  By Proposition \ref{p:open-relation} and \cite{Bourbaki1966}, I.8.3, Proposition 8, $X_W$ is Hausdorff if and only if the graph of $\sim_c$ is closed in $\A_X^W\times\A_X^W$, that is, if whenever $f_n\to f$ and $g_n\to g$ in $\A_X^W$ and $f_n\sim_c g_n$ for all $n$, we have $f\sim_c g$.

Note that $q$ has local continuous sections since $c$ does, and that a $W$-disc structure on $X$ is therefore nothing but a continuous section $\sigma:X\to X_W$ of $\pi$ (so, in particular, $\pi$ is surjective).  Condition (N) says that $\sigma$ extends the tautological section $W\to X_W^0$ that takes a point in $W$ to the class of the constant disc at that point.  Thus we can restate Theorems \ref{t:first-disc-formula} and \ref{t:second-disc-formula} as follows.

\begin{theorem}  \label{t:new-formulation}
Let $W\subset X$ be domains in $\C^n$.  Suppose that one of the following two conditions holds.
\begin{itemize}
\item  There is a continuous section $X\to X_W^0$ of $\pi:X_W\to X$ extending the tauto\-logical section on $W$.
\item  The restriction of $\pi$ to some connected component of $X_W$ is a bijection onto $X$.
\end{itemize}
Then, for every upper semicontinuous function $\phi:W\to[-\infty,\infty)$,
\[ S\phi = E_{\A_X^W}\phi. \]
\end{theorem}

\section{Hartogs domains and Kiselman's minimum principle}

\noindent
Let $Y$ be a domain in $\C^{n-1}$, $n \geq 2$.  Let $r:Y\to[0,\infty)$ be an upper semicontinuous function and $R:Y\to(0,\infty]$ be lower semicontinuous with $r<R$.  Then 
\[ W=\{(z',z_n)\in Y\times\C:r(z')<\lvert z_n\rvert<R(z')\} \] 
is a Hartogs domain in $\C^n$.  It is well known that if $W$ is pseudoconvex, then $\log r$ is plurisubharmonic and $\log R$ is plurisuperharmonic.  Let
\[ X =\{(z',z_n)\in Y\times\C:\lvert z_n\rvert<R(z')\} \]
be the completion of $W$.

Assume $W$ is pseudoconvex.  We claim that every analytic disc $f=(f',f_n)\in\A_X^W$ is centre-homotopic to a vertical disc of a special kind.  Namely, since $f(\T)\subset W$, $f_n$ has no zeros on $\T$.  Let $h$ be a harmonic extension of $\log\lvert f_n\rvert$ to $\overline\D$.  Now $f$ is centre-homotopic in $\A_X^W$ to $\zeta\mapsto f(s\zeta)$ for $s<1$ close enough to $1$, so we may assume that $f$ is smooth on $\overline\D$.  Then $h$ has a harmonic conjugate $k$ on $\D$ that extends continuously to $\overline\D$, so $H=e^{h+ik}:\overline\D\to\C^*$ is continuous with $H|\D$ holomorphic and $\lvert f_n\rvert = \lvert H\rvert$ on $\T$.  Hence $\lvert f_n/H\rvert\leq 1$ on $\overline\D$.  Again since $f(\T)\subset W$, $r\circ f'<\lvert f_n\rvert<R\circ f'$ so $\log r\circ f'<h<\log R\circ f'$ on $\T$.  Since $\log r\circ f'$ is subharmonic and $\log R\circ f'$ superharmonic, $r\circ f'< e^h=\lvert H\rvert<R\circ f'$ on $\overline\D$.

Consider the continuous map $[0,1]\to\A_X$, $t\mapsto f^t$, where
\[ f^t(\zeta) = \big(f'(t\zeta), f_n(\zeta)H(t\zeta)/H(\zeta)\big), \qquad \zeta\in\overline\D. \]
For $\zeta\in\T$, $\lvert f_n^t(\zeta)\rvert = \lvert H(t\zeta)\rvert$, so $f^t(\T)\subset W$ for every $t\in[0,1]$.  Also, $f^t(0)=f(0)$ for all $t\in[0,1]$, and $f^1=f$, so $f$ is centre-homotopic to the vertical disc $g=f^0\in\A_X^W$ with $g(\zeta)=(f'(0),f_n(\zeta)H(0)/H(\zeta))$.  Note that $\lvert g_n\rvert = \lvert H(0)\rvert$ on $\T$.  

There are two cases.  If $f_n$ and therefore $g_n$ has no zeros in $\D$, then $g$ is constant, so in particular $f(0)=g(0)\in W$.  If $f_n$ has $k\geq 1$ zeros in $\D$, then $f_n/H$ is a Blaschke product with $k$ factors.  Hence, either $f(\overline\D)\subset W$ and $f$ is centre-homotopic to a constant disc, or $f$ is centre-homotopic to an analytic disc of the form $(f'(0), sB)\in\A_X^W$, where $B$ is a Blaschke product and $r(f'(0))<s<R(f'(0))$.

This shows that $\A_X^W$ has one connected component $\A_k$ for each winding number $k\geq 0$, and within each component, two discs with the same centre are centre-homotopic.  Clearly, $c(\A_0)=W$ and $c(\A_k) = X$ for $k\geq 1$, where $c$ is the centre map.  The quotient $X_W$ has a component $q(\A_0)$ biholomorphic to $W$, and components $q(\A_k)$, $k\geq 1$, each biholomorphic to $X$.  In particular, $X_W$ is Hausdorff.  Evidently, $X$ is not a schlicht disc extension of $W$.  Finally, Theorem \ref{t:second-disc-formula} implies that our disc formula holds for the pair $W\subset X$.

\smallskip

The disc formula provides a new proof of Kiselman's minimum principle (\cite{Kiselman1978}, Theorem 2.2) in the present setting.  Let $\phi:W\to[-\infty,\infty)$ be plurisubharmonic and rotation-invariant in the last variable.  Define
\[ \psi:Y\to[-\infty,\infty), \quad \psi(z')=\inf\limits_{(z',z_n)\in W}\phi(z',z_n). \]
Let us refer to $\psi$ as the {\it infimum function} of $\phi$.  Kiselman's minimum principle states that $\psi$ is plurisubharmonic.  We will prove it by showing that $\psi(z')=S\phi(z',z_n)$ for all $(z',z_n)\in X\setminus\overline W$.

It is clear that $S\phi(z',z_n)\leq \psi(z')$ for $(z',z_n)\in X\setminus\overline W$.  To prove the opposite inequality, take $z=(z',z_n)\in X$ and $\epsilon>0$, and find $f=(f',f_n)\in\A_X^W$ with $f(0)=z$ and $H_\phi(f)<S\phi(z)+\epsilon$.  As above, there is a continuous function $H:\overline\D\to\C^*$ such that $H|\D$ is holomorphic, $\lvert f_n\rvert=\lvert H\rvert$ on $\T$, and $r\circ f'< \lvert H\rvert<R\circ f'$ on $\overline\D$.  Let
\[ F:\overline\D\times\overline\D\to\C^n, \quad F(\zeta,\xi)=(f'(\zeta), H(\zeta)\xi). \]
Then $F(\overline\D,\T)\subset W$, so for each $\xi\in\T$, $F(\cdot,\xi)$ is an analytic disc in $W$ and
\[ \psi(z')\leq\phi(F(0,\xi)) \leq \int_\T \phi(F(\cdot,\xi))\,d\lambda. \]
If we average over $\xi$, change the order of integration, and use the rotation-invariance of $\phi$, we get
\[ \psi(z') \leq \int_\T\int_\T \phi(F(\zeta,\xi))\,d\lambda(\xi)d\lambda(\zeta) = \int_\T\phi(F(\zeta,f_n(\zeta)/H(\zeta)))\,d\lambda(\zeta)=H_\phi(f).\]
This proves that $\psi(z')\leq S\phi(z',z_n)$.

\smallskip

The observation that Kiselman's infimum function is a plurisubharmonic subextension can be generalised.  Let $W$ be a Hartogs domain in $\C^n=\C^{n-1}\times\C$, such that if $(z_1,\dots,z_n)\in W$ and $\eta\in\T$, then $(z_1,\dots,z_{n-1},\eta z_n)\in W$.  Let $p:W\to\C^{n-1}$ be the projection $(z_1,\dots,z_n)\mapsto (z_1,\dots,z_{n-1})$, and let $Y$ be the domain $p(W)$ in $\C^{n-1}$.  For every $y\in Y$, each vertical fibre $p^{-1}(y)$ is a disjoint union of at most one disc and some number of annuli, possibly none, all centred at the origin in $\{y\}\times\C$.  Define an upper semicontinuous function $r:W\to[0,\infty)$ and a lower semicontinuous function $R:W\to(0,\infty]$ as follows.  If $z\in W$ is contained in an annulus in $p^{-1}(p(z))$, then $r(z)$ is the inner radius and $R(z)$ the outer radius of the annulus.  If $z$ is contained in a disc in $p^{-1}(p(z))$, then $R(z)$ is the radius of the disc and $r(z)=0$.

Suppose $W$ is pseudoconvex.  Then $\log r$ is plurisubharmonic and $\log R$ is plurisuperharmonic.  If $p^{-1}(y_0)$ contains a disc for some $y_0\in Y$, then $p^{-1}(y)$ contains a disc for all nearby $y\in Y$, so $\log r=-\infty$ on an open subset of $W$.  Since $W$ is connected, $\log r=-\infty$ on all of $W$, so every vertical fibre is a disc or a punctured disc.  The latter are ruled out by the Kontinuit\"atssatz.  Thus pseudoconvex Hartogs domains are divided into two classes: the complete ones, whose vertical fibres are discs, and the incomplete ones, whose vertical fibres are unions of annuli.

We assume that $W$ is pseudoconvex and incomplete.  Kiselman showed that identifying each annulus in each vertical fibre of $W$ to a point yields a connected Hausdorff space $A$, which is a pseudoconvex Riemann domain over $Y$ (\cite{Kiselman1978}, Proposition 2.1 and Corollary 2.3).  Now $r$ and $R$ induce functions on $A$, plurisubharmonic and plurisuperharmonic, respectively, and $W$ is reincarnated as the Hartogs domain
\[ W=\{(a,\zeta)\in A\times\C:r(a)<\lvert\zeta\rvert<R(a)\} \]
over $A$, with completion
\[ X=\{(a,\zeta)\in A\times\C:\lvert\zeta\rvert<R(a)\}, \]
in which $A$ is embedded as $A\times\{0\}$.  

Let $\phi:W\to[-\infty,\infty)$ be plurisubharmonic and rotation-invariant in the sense that $\phi(a,\eta\zeta)=\phi(a,\zeta)$ for all $(a,\zeta)\in W$ and $\eta\in\T$.  Define the infimum function
\[ \psi:A\to[-\infty,\infty), \quad \psi(a)=\inf\limits_{(a,\zeta)\in W}\phi(a,\zeta).\]  
Kiselman's minimum principle states that $\psi$ is plurisubharmonic.  Our generalisation of the principle is as follows.

\begin{theorem}  \label{t:generalisation}
Let $W$ be an incomplete pseudoconvex Hartogs domain in $\C^n$ with completion $X$ over the Kiselman quotient $A$ of $W$.
\begin{enumerate}
\item  If $\phi$ is a plurisubharmonic and rotation-invariant function on $W$, then the infimum function $\psi$ of $\phi$ is the restriction to $A$ of the envelope $E_{\A_X^W}\phi$.
\item  For every upper semicontinuous function $\phi$ on $W$, the envelope $E_{\A_X^W}\phi$ is pluri\-subharmonic on $X$, and is therefore the largest plurisubharmonic subextension of $\phi$ to $X$.
\end{enumerate}
\end{theorem}

We note that we are not giving a new proof of Kiselman's minimum principle, because the principle is used in the proof that $A$ is pseudoconvex, which implies that $X$ is Stein.  We need this in our proof of (2) in order to apply Theorem \ref{t:second-disc-formula}.  Statement (1), proved from scratch by the method used above, shows that (2) does indeed generalise Kiselman's minimum principle.

\begin{proof}
(1)  For $a\in A$, $\psi(a)$ is the infimum of $H_\phi(f)$ over all analytic discs $f\in\A_X^W$ of the form $f(\zeta)=(a,s\zeta)$ with $r(a)<s<R(a)$.  Hence $\psi\geq E_{\A_X^W}\phi|A$.  The opposite inequality can be proved as above, assuming, as we may, that $f$ is smooth on $\overline\D$.

(2)  Let $\B$ be the connected component of $\A_X^W$ containing an analytic disc of the form $\zeta\mapsto (a,s\zeta)$ with $a\in A$ and $r(a)<s<R(a)$.  Then $\B$ contains all analytic discs of the form $\zeta\mapsto (a,s\alpha(\zeta))$ with $a\in A$, $r(a)<s<R(a)$, and $\alpha\in\operatorname{Aut}\D$, so $\B$ covers $X$.  As above, we can prove that every analytic disc in $\B$ is centre-homotopic to an analytic disc of the latter form.  Therefore, if two analytic discs in $\B$ have the same centre, then they are centre-homotopic.  Since $X$ is Stein, Theorem \ref{t:second-disc-formula} now implies that our disc formula holds for the pair $W\subset X$.
\end{proof}


\newpage

\begin{thebibliography}{99}

\bibitem{Bourbaki1966}
N. Bourbaki.  {\it General topology.}  Hermann, Addison-Wesley, 1966.

\bibitem{BuSchachermayer1992} 
Bu, S.\ Q.\ and W.\ Schachermayer.  {\it Approximation of Jensen measures by image measures under holomorphic functions and applications.}  Trans.\ Amer.\ Math.\ Soc.\ \textbf{331} (1992) 585--608.

\bibitem{DrnovsekForstneric2011a}
Drinovec Drnov\v sek, B.\ and F.\ Forstneri\v c.  {\it The Poletsky-Rosay theorem on singular complex spaces,} preprint, 2011, arXiv 1104.3968.

\bibitem{DrnovsekForstneric2011b}
Drinovec Drnov\v sek, B.\ and F.\ Forstneri\v c.  {\it Disc functionals and Siciak-Zaharyuta extremal functions on singular varieties,} preprint, 2011, arXiv 1109.3947.

\bibitem{Edigarian1997}
Edigarian, A.  {\it On definitions of the pluricomplex Green function.}  Ann.\ Polon.\ Math.\ \textbf{67} (1997) 233--246.

\bibitem{Forstneric2007}
Forstneri\v c, F.  {\it Manifolds of holomorphic mappings from strongly pseudoconvex domains.}  Asian J.\ Math.\ \textbf{11} (2007) 113--126.

\bibitem{Goluzin1969}
Goluzin, G.\ M.  {\it Geometric theory of functions of a complex variable.}  Translations of Mathematical Monographs, 26.  Amer.\ Math.\ Soc., 1969.

\bibitem{Kiselman1978}
Kiselman, C.\ O.  {\it The partial Legendre transformation for plurisubharmonic functions.}  Invent.\ Math.\ \textbf{49} (1978) 137--148.

\bibitem{LarussonSigurdsson1998}
L\'arusson, F.\ and R.\ Sigurdsson.  {\it Plurisubharmonic functions and analytic discs on manifolds.}  J.\ reine angew.\ Math.\ \textbf{501} (1998) 1--39.

\bibitem{LarussonSigurdsson2003}
L\'arusson, F.\ and R.\ Sigurdsson.  {\it Plurisubharmonicity of envelopes of disc functionals on manifolds.}  J.\ reine angew.\ Math.\ \textbf{555} (2003) 27--38.

\bibitem{LarussonSigurdsson2005}
L\'arusson, F.\ and R. Sigurdsson.  {\it The Siciak-Zahariuta extremal function as the envelope of disc functionals.}  Ann.\ Polon.\ Math.\ \textbf{86} (2005) 177--192.

\bibitem{LarussonSigurdsson2009}
L\'arusson, F.\ and R.\ Sigurdsson.  {\it Siciak-Zahariuta extremal functions, analytic discs and polynomial hulls.}  Math.\ Ann.\ \textbf{345} (2009) 159--174.

\bibitem{Magnusson2011}
Magn\'usson, B.\ S.  {\it Extremal $\omega$-plurisubharmonic functions as envelopes of disc functionals.}  Arkiv Math.\ \textbf{49} (2011) 383--399.

\bibitem{Poletsky1991}
Poletsky, E.\ A.  {\it  Plurisubharmonic functions as solutions of variational problems.}  Several complex variables and complex geometry (Santa Cruz CA, 1989), 163--171.  Proc.\ Sympos.\ Pure Math., 52, Part 1.  Amer.\ Math.\ Soc., 1991.

\bibitem{Poletsky1993}
Poletsky, E.\ A.  {\it Holomorphic currents.}  Indiana Univ.\ Math.\ J.\ \textbf{42} (1993) 85--144.

\bibitem{PoletskyShabat1989}
Poletsky, E.\ A.\ and B.\ V.\ Shabat.  {\it Invariant metrics.}  Several Complex Variables III, 63--111.  Encyclopaedia of Mathematical Sciences, 9.  Springer-Verlag, 1989.

\bibitem{Rado1922}
Rad\'o, T.  {\it Sur la repr\'esentation conforme de domaines variables.}  Acta litterarum ac scientiarum regiae universitatis hungaricae Francisco-Josephinae. Sectio scientiarum mathematicarum \textbf{1}  (1922--1923) 180--186.

\bibitem{RashkovskiiSigurdsson2005}
Rashkovskii, A.\ and R.\ Sigurdsson.  {\it Green functions with singularities along complex spaces.}  Internat.\ J.\ Math.\ \textbf{16} (2005) 333--355.

\bibitem{Rosay2003}
Rosay, J.-P.  {\it Poletsky theory of disks on holomorphic manifolds.}
 Indiana Univ.\ Math.\ J.\ \textbf{52} (2003) 157--169.

\end{thebibliography}
\end{document}